\documentclass{amsart}
\usepackage{amsmath,amssymb}
\usepackage{graphicx}
\usepackage{tikz}
\usetikzlibrary{arrows.meta,decorations.markings}
\usepackage{microtype}
\usepackage{xcolor}
\definecolor{figureteal}{HTML}{287D8E}
\definecolor{figuregold}{HTML}{C48736}
\definecolor{figureviolet}{HTML}{7B6AA6}
\definecolor{figureink}{HTML}{405466}
\usepackage[hidelinks]{hyperref}
\DeclareTextCommand{\k}{OT1}[1]{\UseTextAccent{T1}{\k}{#1}}

\allowdisplaybreaks[4]

\newtheorem{thm}{Theorem}
\newtheorem{lem}[thm]{Lemma}

\newtheorem{addendum}[thm]{Addendum}
\theoremstyle{definition}
\newtheorem{rmk}[thm]{Remark}
\newtheorem{question}[thm]{Question}

\newcommand{\Z}{\mathbb{Z}}
\newcommand{\R}{\mathbb{R}}
\newcommand{\x}{\times}

\newcommand{\Int}{\mathrm{int}}

\begin{document}

\title[Exotic smooth structures on nonorientable 4-manifolds]
{Exotic smooth structures on nonorientable 4-manifolds}
\author[R. \.{I}. Baykur]{R. \.{I}nan\c{c} Baykur}
\address{Department of Mathematics and Statistics, University of Massachusetts, Amherst, MA 01003, USA}
\email{inanc.baykur@umass.edu}
\date{}

\begin{abstract}
We construct closed nonorientable aspherical $4$--manifolds with an
arbitrarily large number of distinct smooth structures. This uses
the reflection trick of Davis--Hayden--Huang--Ruberman--Sunukjian.
Our examples become mutually diffeomorphic after either type of
stabilization: taking a product with $\R$, or taking a connected sum
with sufficiently many copies of $S^2\x S^2$.
\end{abstract}

\maketitle
\vspace{-2\baselineskip}
\setcounter{secnumdepth}{2}

\section{Introduction}

Nonorientable exotica in dimension four have an even longer history
than their orientable counterparts, going back to Kreck's work
in 1982~\cite{Kreck1984} and Cappell--Shaneson's earlier homotopy
$\mathbb{RP}^4$'s~\cite{CappellShaneson1976}, whose homeomorphism
to $\mathbb{RP}^4$ was established later~\cite[p.~221]{Ruberman1984}.
Many interesting constructions, invariants, and classification results
have appeared over the years, in work of
Fintushel--Stern~\cite{FintushelStern1981,FintushelStern1984},
Akbulut~\cite{Akbulut1984,Akbulut1988},
Gompf~\cite{Gompf1984,Gompf1991}, Ruberman~\cite{Ruberman1984}, Gilkey~\cite{Gilkey1985},
Stolz~\cite{Stolz1988}, Oledzki~\cite{Oledzki1991},
Hambleton--Kreck--Teichner~\cite{HambletonKreckTeichner1994},
Wang~\cite{Wang1995},
Torres~\cite{Torres2017,TorresPrescribed2026,TorresClassification2026,Torres2026},
Bais--Torres~\cite{BaisTorres2025,BaisTorresCutPaste2025}, Debray~\cite{Debray2022}, and
Kasprowski--Powell~\cite{KasprowskiPowell2025}.
Yet, surprisingly, while orientable exotica often come in
infinite families~\cite{FintushelStern1998,FintushelStern2009,BaykurStipsiczSzabo2024},
to the best of our knowledge no closed nonorientable $4$--manifold
was previously known to admit more than two distinct smooth
structures; cf.~\cite{Torres2026}.

One might hope to obtain nonorientable examples by taking free
quotients of orientable $4$--manifolds with many distinct smooth
structures. The orientable manifolds must then have signature zero.
Such families are in fact available in the literature; see e.g.~\cite{Park2002,AkhmedovParkSakalli2023,BaykurHamada2024,BaykurHamada2026,BaykurStipsiczSzabo2024},
but finding variations of these examples that admit free orientation
reversing involutions has proven challenging---at least to the author
of this article. We thus try a different source for this equivariant
approach, following the recent work on exotic orientable aspherical
$4$--manifolds by Davis--Hayden--Huang--Ruberman--Sunukjian
~\cite{DHHRS-2026}. We obtain arbitrarily
many nonorientable exotica, with asphericity (and thus irreducibility) as an added feature:

\begin{thm}\label{thm:main}
For every positive integer $N$, there exists a closed nonorientable
aspherical topological $4$--manifold $Z_N$ which admits at least $N$
pairwise nondiffeomorphic smooth structures.
These become mutually diffeomorphic after taking a product with
$\R$, and hence after taking a connected sum with sufficiently
many copies of $S^2\x S^2$.
\end{thm}

Here is an outline of both the proof and our paper.
The reflection construction in~\cite{DHHRS-2026} uses
\emph{Davis's reflection group trick}~\cite{Davis1983} with compact pieces
adapted from Hayden--Piccirillo~\cite{HaydenPiccirillo2025} as
\emph{chambers}.
In Section~\ref{sec:construction}, we switch ingredients and choose
compact pieces $X_i$ from Akbulut--Yasui~\cite{AkbulutYasui2013}.
They are pairwise homeomorphic, are homotopy equivalent to $T^2$,
have strictly increasing minimum genera among nonzero homology classes,
and embed smoothly in $D^4$---properties all essential to our proof.
When we use these chambers for an aspherical tiling, the minimal
genus distinction must also hold for surfaces that pass through many
chambers. Take a homologically nonzero surface in a finite boundary
connected sum and choose a summand in which its class has a nonzero
component. The ball embeddings let us embed the whole boundary sum
into this summand, with the other chambers inside disjoint $4$--balls.
Since these balls have zero second homology, the image surface
represents just the nonzero class coming from the chosen summand.
That chamber's genus lower bound therefore applies to the original surface.
This works even when the chambers occur with both orientations, as
in~\cite{DHHRS-2026}.
The reflection construction then assembles copies of these pieces
into noncompact reflection manifolds $D_i$ and chooses closed
nonorientable quotients $Y_i$. Since $X_i\simeq T^2$, the chambers
are aspherical; asphericity then passes through the finite boundary
sums to $D_i$ and $Y_i$. The homeomorphisms between the
chambers extend to the $D_i$ and descend to the $Y_i$.
The characteristic subgroup lemma~\cite[Lemma~5.3]{DHHRS-2026}
then forces hypothetical diffeomorphisms between the $Y_i$ to lift
to the $D_i$, where the genus bounds rule them out.
We end in Section~\ref{sec:remarks} with the still open
Question~\ref{qst:infinite} on
infinitely many smooth structures on a single closed nonorientable
$4$--manifold, which originally motivated this work.
In Remark~\ref{rmk:earlier-methods}, we share our notes on why earlier
constructions and invariants have not already produced arbitrarily
many nonorientable exotica.

\newpage
\subsection*{Conventions and preliminary definitions}
We work with integral homology coefficients. We write $\Z_n$ for the
cyclic group of order $n\ge1$, $\natural$ for boundary connected sum,
and $-X$ for an oriented manifold $X$ with the opposite orientation.
For manifolds, $\cong$, $\approx$ and $\simeq$ denote diffeomorphism,
homeomorphism and homotopy equivalence, respectively. For oriented
manifolds, the first two denote orientation preserving comparisons.
PL abbreviates piecewise linear, and DHHRS refers to the authors
of~\cite{DHHRS-2026}. For a smooth $4$--manifold $X$, the \emph{minimal genus function}
$g_X$ assigns to $u\in H_2(X)$ the minimum genus of a smoothly
embedded, closed, connected, oriented surface in $\Int X$
representing $u$.

When $H_2(X)\ne0$, we also use
\[
\mu(X)=\min\{g_X(u)\mid u\in H_2(X),\ u\ne0\}.
\]
Note that $\mu(-X)=\mu(X)$.
In our notation, DHHRS distinguish their noncompact reflection
manifolds by comparing $\mu$~\cite[Proof of Theorem~4.1]{DHHRS-2026}.
One of their aspherical manifold pairs contains a homologically nonzero torus, whereas every homologically
nonzero surface in the other has genus at least $2$.

Smooth manifolds $X$ and $Y$ are \emph{exotic} copies of one another
if $X\approx Y$ but $X\not\cong Y$.
A manifold is \emph{aspherical} if its universal cover is contractible.
A triangulation is \emph{flag} if every set of pairwise adjacent
vertices spans a simplex, and \emph{flag-no-square} if it is flag
and every $4$-cycle in its $1$-skeleton has at least one of its two
diagonals. For background on the compact Stein surfaces used as our
chambers, the reader can turn to~\cite{OzbagciStipsicz2004}; for
Coxeter groups and the reflection group trick, to~\cite{Davis1983};
and for flag-no-square triangulations, to~\cite{PrzytyckiSwiatkowski2009}.
We follow the reflection construction as presented in~\cite{DHHRS-2026}.

\section{The construction}\label{sec:construction}

We begin by choosing compact pieces whose genus bounds survive the
reflection trick.

\begin{lem}\label{lem:chambers}
For every positive integer $N$, there are compact connected oriented smooth
$4$--manifolds $X_1,\ldots,X_N$ with the following properties:
\begin{enumerate}
\item $X_i\approx X_j$ rel a fixed collar of their identified boundaries;
\item $X_i\simeq T^2$ and $X_i$ embeds smoothly in $D^4$;
\item $\mu(X_1)<\cdots<\mu(X_N)$.
\end{enumerate}
\end{lem}

\begin{proof}
Apply~\cite[Theorem~5.16]{AkbulutYasui2013} to $X=T^2\x D^2$,
taking $N$ Stein members with the parameter choices described below.
Part~(5) of that theorem embeds each $X_i$ smoothly in $X$, which in this case
itself embeds in $D^4$ as a tubular neighborhood of the standard
unknotted torus. Part~(6) realizes each $X_i$ by twisting a cork
$C_i\subset\Int X_0$ in a common manifold $X_0$.
The resulting homeomorphisms $X_0\to X_i$ are the identity outside
the corks, hence on a common boundary collar. This proves the first claim.

To check the homotopy type, note that $X_0$ is obtained from
$X=T^2\x D^2$ by $W^{-}$-modifications, each of which attaches
a homotopically canceling pair of $1$- and $2$-handles.
Thus $X_0\simeq X$, and since each $X_i$ is homeomorphic to $X_0$,
we obtain $X_i\simeq T^2$; see~\cite[Theorem~3.7]{YasuiGenus2025}.

For the genus bounds, Akbulut--Yasui's calculations
in~\cite[Lemmas~5.13 and~5.14]{AkbulutYasui2013}
start with fixed integers $g_0>0$ and $c$, which depend on the initial Legendrian
handle presentation of $T^2\x D^2$, and integers
 $0<p_1<\cdots<p_N$,
with generators $u_i\in H_2(X_i)\cong\Z$ and Stein structures $J_i$
satisfying
\[
u_i^2=0,\qquad
g_{X_i}(u_i)\le g_0+p_i,\qquad
|\langle c_1(J_i),u_i\rangle|=2p_i+c.
\]

Take $p_1$ sufficiently large for the initial requirements of the
construction in~\cite{AkbulutYasui2013} and so that $2p_1+c>0$. For every $k\ne0$, the
adjunction inequality for Stein surfaces~\cite{OzbagciStipsicz2004}
then gives
\[
2g_{X_i}(k u_i)-2
\ge (k u_i)^2+|\langle c_1(J_i),k u_i\rangle|
= |k|(2p_i+c)\ge 2p_i+c.
\]
Together with the upper bound, this yields
\begin{equation}\label{eq:genus-intervals}
p_i+1+\frac{c}{2}\le\mu(X_i)\le g_0+p_i.
\end{equation}
Each $p_i$ can be taken arbitrarily large after the preceding parameters
have been fixed, so we can make these intervals ordered and disjoint.
Thus $\mu(X_1)<\cdots<\mu(X_N)$.\qedhere
\end{proof}

We now apply the reflection construction of DHHRS to these
chambers. We recall its relevant features, then choose free orientation
reversing involutions on the orientable quotients.

Fix $N$ and the manifolds of Lemma~\ref{lem:chambers}. Choose a smooth
triangulation of their common boundary and subdivide it to obtain a
flag-no-square triangulation $\mathcal T$; such a subdivision exists
by~\cite[Proposition~2.13]{PrzytyckiSwiatkowski2009}.
Following~\cite[Section~3.1]{DHHRS-2026}, write
$V$ for its vertex set and let
\[
W=\langle s_v\ \text{for }v\in V \mid s_v^2=1,\quad
s_vs_w=s_ws_v\text{ if }\{v,w\}\text{ is an edge in }\mathcal T\rangle
\]
be the associated right-angled Coxeter group.

\begin{figure}[htb]
\centering
\begin{tikzpicture}[font=\small, line join=round, draw=figureink,
  tauarrow/.style={Stealth-Stealth, semithick, draw=figureviolet}]
\node at (2.3,4.1) {$D_i/C$};
\node at (9.5,4.1) {$Y_i=D_i/H$};

\draw[fill=figureteal!18] (0,2.4) rectangle (2.3,3.5);
\draw[fill=figuregold!20] (2.3,2.4) rectangle (4.6,3.5);
\draw[fill=figuregold!20] (0,0.5) rectangle (2.3,1.6);
\draw[fill=figureteal!18] (2.3,0.5) rectangle (4.6,1.6);
\draw[very thick] (2.3,2.4)--(2.3,3.5);
\draw[very thick] (2.3,0.5)--(2.3,1.6);
\node[align=center] at (1.15,2.95) {$0$\\[2pt]$+X_i$};
\node[align=center] at (3.45,2.95) {$e_v$\\[2pt]$-X_i$};
\node[align=center] at (1.15,1.05) {$a$\\[2pt]$-X_i$};
\node[align=center] at (3.45,1.05) {$a+e_v$\\[2pt]$+X_i$};
\node at (2.3,3.72) {$v$};
\node at (2.3,1.82) {$v$};
\draw[tauarrow] (1.15,1.72)--(1.15,2.28)
  node[midway,left=3pt] {$\tau$};
\draw[tauarrow] (3.45,1.72)--(3.45,2.28)
  node[midway,right=3pt] {$\tau$};
\node at (2.3,0.02) {chambers indexed by $A$};

\draw[-{Stealth},thick,draw=figureviolet] (5.02,2)--(6.78,2)
  node[midway,above=4pt] {$q$};
\node at (5.9,1.52) {$b\sim b+a$};

\draw[fill=figureviolet!14] (7.2,1.45) rectangle (9.5,2.55);
\draw[fill=figureviolet!14] (9.5,1.45) rectangle (11.8,2.55);
\draw[very thick] (9.5,1.45)--(9.5,2.55);
\node[align=center] at (8.35,2) {$[0]$\\[2pt]$X_i$};
\node[align=center] at (10.65,2) {$[e_v]$\\[2pt]$X_i$};
\node at (9.5,2.77) {$v$};
\node at (8.35,0.92) {$0\sim a$};
\node at (10.65,0.92) {$e_v\sim a+e_v$};
\node at (9.5,0.02) {chambers indexed by $A/\langle a\rangle$};

\end{tikzpicture}
\caption{Schematic chamber identifications. In $D_i/C$, the chambers
are indexed by $b\in A$, and the deck involution $\tau$ sends the
chamber indexed by $b$ to that indexed by $b+a$. In $Y_i=D_i/H$,
they are indexed by the cosets $[b]\in A/\langle a\rangle$.
Only adjacency across a $v$-mirror is shown, where $e_v=\rho(s_v)$.}
\label{fig:involution}
\end{figure}
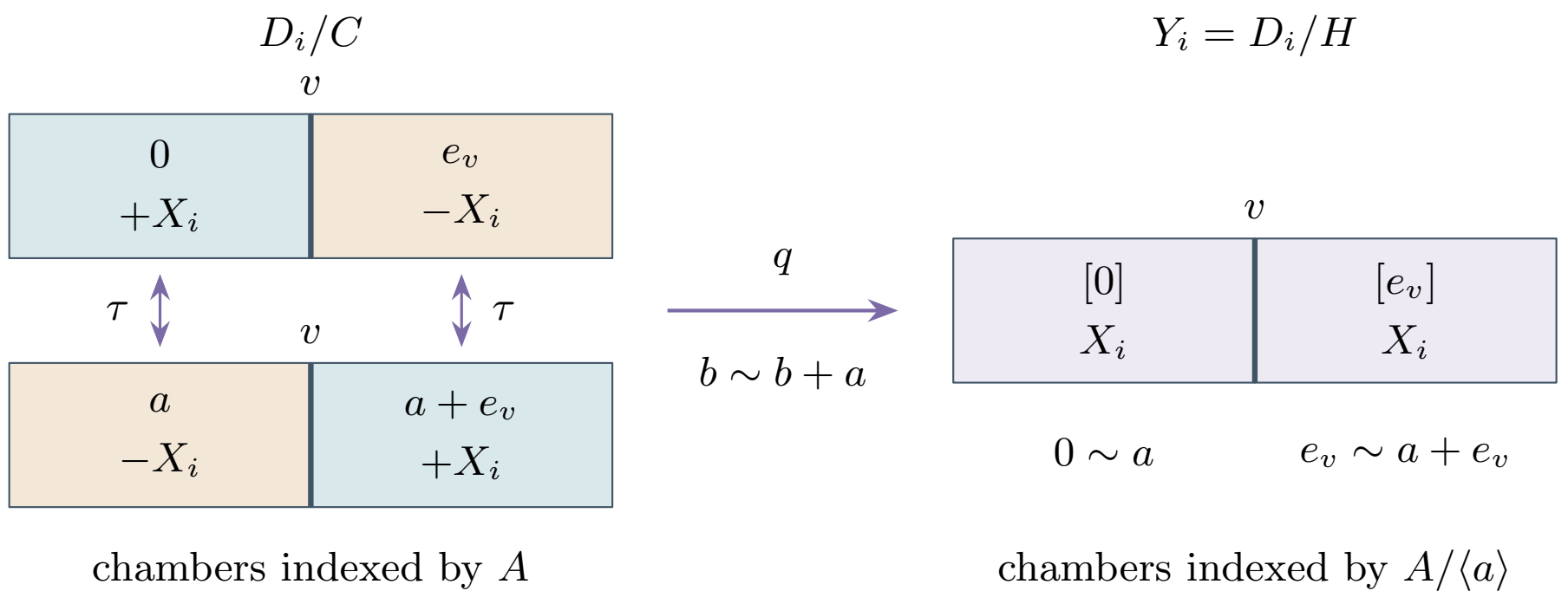

For each $v\in V$, the closed $3$--cell $M_{i,v}\subset\partial X_i$
dual to $v$ in this triangulation is called a \emph{mirror}.
Following Davis~\cite{Davis1983}, the reflection manifold $D_i=D(X_i)$
is built from copies $wX_i$ indexed by $w\in W$, identifying the
$v$-mirror in $wX_i$ with that in $ws_vX_i$ by the identity map. We give adjacent chambers opposite orientations, so the gluing
reverses their induced boundary orientations. Note that the first relation in $W$ records that reflecting twice returns to the same
chamber, whereas the second allows the two reflections to commute where the
corresponding mirrors meet. See the left side of Figure~\ref{fig:involution} for a schematic
of adjacency across a mirror in the quotient $D_i/C$ introduced below.

More precisely, for $x\in X_i$ put
$W_x=\langle s_v\mid x\in M_{i,v}\rangle \subset W$, with $W_x=\{1\}$
when $x\in\Int X_i$. Then
\[
D_i=(W\times X_i)/\sim,\qquad
(w,x)\sim(w',x')\ \Longleftrightarrow
x=x'\ \text{and}\ w^{-1}w'\in W_x.
\]
Writing $wx=[w,x]$, left multiplication gives the action
$h[w,x]=[hw,x]$. Keeping the coordinate $x$ fixed means that we move
to the corresponding point in the chamber $hwX_i$.
This action is proper and smooth~\cite[Section~3.2]{DHHRS-2026}.
The stabilizer of $[w,x]$ is $wW_xw^{-1}$. The mirrors through $x$
correspond to vertices of a simplex of $\mathcal T$, so $W_x\cong(\Z_2)^r$
with $r\le4$. (Four is the maximum possible, realized by the mirrors dual to the four vertices
of a tetrahedron, which meet at its barycenter.)
The chamber $wX_i$ has orientation $(-1)^{\ell(w)}$ times that of
$X_i$, where $\ell(w)$ is the word length, i.e. the least number of generators needed to
express $w$.

Ordering chambers by nondecreasing word length gives a compact
exhaustion $P_{i,n}$ of $D_i$, where each $P_{i,n}$ is the union of the first
$n$ chambers and $D_i=\bigcup_n\Int P_{i,n}$.
By~\cite[Lemmas~3.1--3.3]{DHHRS-2026}, after rounding corners,
\begin{equation}\label{eq:boundary-sum}
P_{i,n}\approx_{\mathrm{PL}} X_i^{(1)}\natural\cdots\natural X_i^{(n)},
\end{equation}
where each summand is a copy of $X_i$ with one of the two orientations.
The chamber and exhaustion inclusions induce injections on $H_2$.
Indeed, attaching each chamber along a $3$--ball adds its $H_2$
as a direct summand by Mayer--Vietoris. These inclusions remain injective
in $D_i$, since every bounding chain lies in some $P_{i,n}$.
These finite boundary sums are homotopy equivalent to wedges of
aspherical chambers. Asphericity therefore passes to $D_i$, since
every map $S^k\to D_i$ with $k\ge2$ has image in some $P_{i,n}$.

\begin{figure}[htb]
\centering
\begin{tikzpicture}[x=0.92cm,y=0.75cm,font=\small,
  chamber/.style={draw=figureink, thick, fill=white},
  handle/.style={draw=figureink!80, line width=0.5pt, double=white, double distance=3pt}]
\node at (2.2,2.25) {$X_i^{(j)}\natural X_i^{(k)}\natural X_i^{(\ell)}$};
\node at (10,2.25) {$X_i^{(k)}$};

\draw[handle] (0.5,0)--(1.25,0);
\draw[handle] (3.15,0)--(3.9,0);
\draw[chamber,draw=figuregold!85!black,fill=figuregold!18] (0,0) ellipse (0.68 and 0.57);
\draw[chamber,draw=figureteal,fill=figureteal!18] (2.2,0) ellipse (1.03 and 0.83);
\draw[chamber,draw=figureviolet,fill=figureviolet!16] (4.4,0) ellipse (0.68 and 0.57);
\node at (0,0) {$X_i^{(j)}$};
\node at (2.2,0) {$X_i^{(k)}$};
\node at (4.4,0) {$X_i^{(\ell)}$};

\node[text=figureink] at (6.025,0.2) {\Large$\lhook\joinrel\longrightarrow$};

\draw[thick,draw=figureteal!85!black,fill=figureteal!3] (10,0) ellipse (3.05 and 1.9);
\draw[dashed,draw=figureink!75] (8.88,-1.28) circle (0.42);
\draw[dashed,draw=figureink!75] (11.12,-1.28) circle (0.42);
\draw[handle] (8.88,-0.65)--(8.88,-1.10);
\draw[handle] (11.12,-0.65)--(11.12,-1.10);
\draw[chamber,draw=figureteal,fill=figureteal!18] (10,0) ellipse (2.1 and 0.88);
\node at (10,0) {copy of $X_i^{(k)}$};
\draw[chamber,draw=figuregold!85!black,fill=figuregold!18] (8.88,-1.28) ellipse (0.36 and 0.31);
\draw[chamber,draw=figureviolet,fill=figureviolet!16] (11.12,-1.28) ellipse (0.36 and 0.31);
\node[font=\fontsize{6}{7}\selectfont,inner sep=0pt] at (8.88,-1.28) {$X_i^{(j)}$};
\node[font=\fontsize{6}{7}\selectfont,inner sep=0pt] at (11.12,-1.28) {$X_i^{(\ell)}$};
\end{tikzpicture}
\caption{Embedding a boundary connected sum into the chosen chamber
$X_i^{(k)}$, illustrated for three summands. The summands $X_i^{(j)}$
and $X_i^{(\ell)}$ lie in disjoint $4$--balls (dashed circles)
in its boundary collar. }
\label{fig:genus-embedding}
\end{figure}
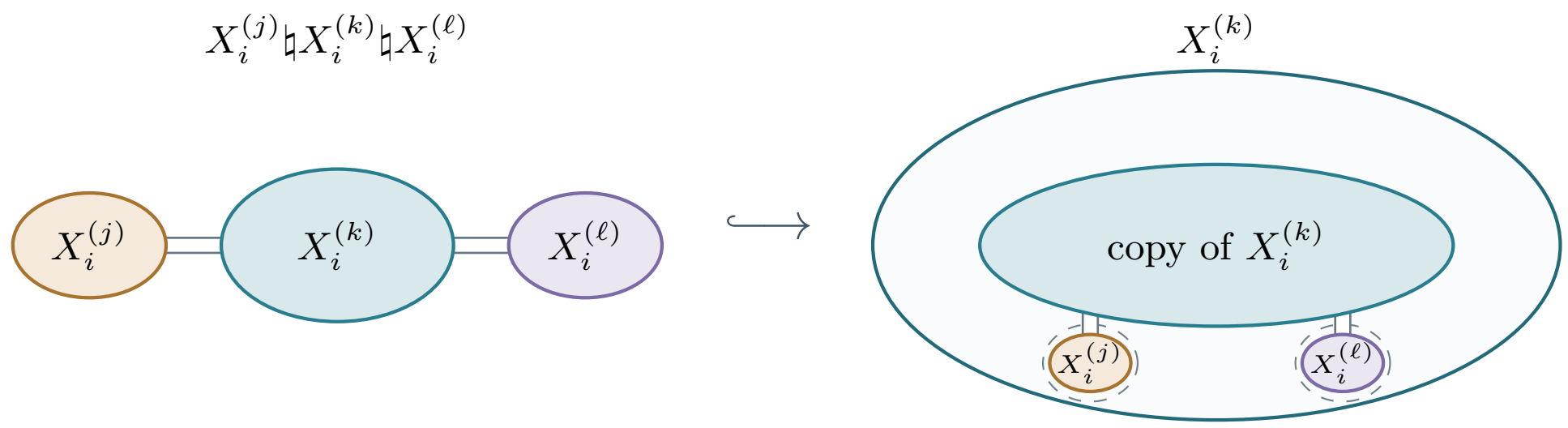

The genus comparison in the proof of~\cite[Theorem~4.1]{DHHRS-2026}
gives, for our chambers,
\begin{equation}\label{eq:genus-reflection}
\mu(D_i)=\mu(X_i).
\end{equation}
Here the chamber inclusion gives $\mu(D_i)\le\mu(X_i)$.
For the reverse inequality, a homologically nonzero smooth surface
in $D_i$ lies in the interior of some $P_{i,n}$. The PL identification
above, followed by smoothing~\cite[Theorem~2]{Wall1967} as in that
proof, gives a smooth surface of the same genus in the boundary sum.
Choose a summand $X_i^{(k)}$ in which its homology class has a nonzero
component. The ball embeddings of Lemma~\ref{lem:chambers} let us
embed the whole boundary sum in $X_i^{(k)}$, with the other summands
in disjoint $4$--balls in a boundary collar; reflect the balls when
needed to match orientations. This embedding induces projection
onto $H_2(X_i^{(k)})$, so the image surface is homologically nonzero
and has genus at least $\mu(X_i^{(k)})=\mu(X_i)$.
Figure~\ref{fig:genus-embedding} illustrates this.

We next obtain nonorientable examples by finding free orientation
reversing involutions on closed orientable quotients of the $D_i$.
Consider the abelianization map
\[
\rho\colon W\longrightarrow A=\bigoplus_{v\in V}\Z_2,
\qquad \rho(s_v)=e_v,\qquad C=\ker\rho=[W,W],
\]
where $e_v$ is the basis vector corresponding to $v$. The $v$-coordinate
of $\rho(w)$ records the parity of the number of occurrences of $s_v$
in a word representing $w$.

\begin{lem}\label{lem:nonorientable-quotients}
The manifolds $D_i/C$ are closed, orientable and aspherical, and
admit free orientation reversing involutions $\tau_i$. Their
quotients $Y_i=(D_i/C)/\langle\tau_i\rangle$ have the form $D_i/H$
for a common torsion-free subgroup $H$ of finite index in $W$.
\end{lem}

\begin{proof}

By~\cite[Theorem~12.3.4(i)]{Davis2025}, every finite subgroup
of $W$ is conjugate into a finite subgroup generated by some of the
$s_v$. These generators commute, since two noncommuting generators
generate an infinite dihedral group, and their vertices span a simplex
because $\mathcal T$ is flag. Hence every nontrivial finite-order
element of $W$ is conjugate to a product of distinct generators
from a simplex. Its image under $\rho$ is therefore nonzero and has
at most four nonzero coordinates (since $\dim\mathcal T=3$). In particular, $C$ is torsion-free~\cite{DHHRS-2026}.
In a word representing an element of $C$, each generator
occurs an even number of times. Since each generator reverses
orientation, $C$ preserves orientation. Thus $D_i/C$ is a closed
orientable aspherical manifold.

Under the induced $A$-action on $D_i/C$, the stabilizer of
a $C$-orbit $[x]$ is the $\rho$-image of the stabilizer of $x$.
Since point stabilizers in $W$ are finite, the first observation
shows that every element of this image has at most four nonzero coordinates.

Choose five distinct vertices $v_1,\ldots,v_5\in V$ and put
\[
a=e_{v_1}+\cdots+e_{v_5}.
\]
Since $a$ has five nonzero coordinates, it fixes no point of $D_i/C$.
Also, $2a=0$ in $A$, so the action of $a$ on $D_i/C$ is a free
involution $\tau_i$.
This involution is induced by $s_{v_1}\cdots s_{v_5}$, a composition
of five reflections. Each reflection reverses orientation, so
$\tau_i$ does too.

Put $H=\rho^{-1}(\langle a\rangle)$, which has index $2^{|V|-1}$
in $W$. By the first observation, no nontrivial finite-order element
can have image $0$ or $a$, so $H$ is torsion-free. We obtain closed
nonorientable aspherical manifolds
\[
Y_i=(D_i/C)/\langle\tau_i\rangle=D_i/H,
\]
with orientation double covers $D_i/C$. The chamber pairing is
shown in Figure~\ref{fig:involution}.\qedhere

\end{proof}

\begin{proof}[Proof of Theorem~\ref{thm:main}]
Take $Y_i=D_i/H$ from Lemma~\ref{lem:nonorientable-quotients}.

The homeomorphisms $h_{ij}\colon X_i\to X_j$ of Lemma~\ref{lem:chambers}
are fixed near the identified boundaries. Thus the chamberwise formula
$F_{ij}(wx)=wh_{ij}(x)$ defines $W$-equivariant homeomorphisms
$D_i\to D_j$, which descend to homeomorphisms $Y_i\to Y_j$,
as in the proof of~\cite[Theorem~4.1]{DHHRS-2026}.

The covering subgroup $K_i=\pi_1(D_i)\subset\pi_1(Y_i)$ is
characteristic by~\cite[Lemma~5.3]{DHHRS-2026}, since
$\pi_1(X_i)\cong\Z^2$, $\mathcal T$ is flag-no-square, and $H$
is torsion-free of finite index in $W$.

Suppose $g\colon Y_i\to Y_j$ is a diffeomorphism, and let
$f\colon Y_i\to Y_j$ be one of the homeomorphisms above.
The automorphism $(f^{-1}\circ g)_*$ preserves $K_i$, while
$f_*(K_i)=K_j$ by construction. Hence $g_*(K_i)=K_j$, so $g$
lifts to a diffeomorphism $D_i\to D_j$, as in the proof
of~\cite[Theorem~1.1]{DHHRS-2026}.

By~\eqref{eq:genus-reflection}, this would give
\[
\mu(X_i)=\mu(D_i)=\mu(D_j)=\mu(X_j),
\]
which is impossible for $i\ne j$ by Lemma~\ref{lem:chambers}.

The homeomorphisms $Y_i\to Y_j$ above are smooth outside finitely
many contractible corks. Their Casson--Sullivan invariants vanish since
$H^3(C,\partial C;\Z_2)=0$ for each cork $C$, so the $Y_i\x\R$
are mutually diffeomorphic~\cite[Proposition~5.4]{HausmannJahren2018}.
By Gompf~\cite{Gompf1984}, the oriented chambers become diffeomorphic
relative to their boundaries after taking connected sums with copies
of $S^2\x S^2$. These diffeomorphisms agree on the mirrors and glue to stable
diffeomorphisms of the closed quotients.
Taking the common homeomorphism type of the $Y_i$ as $Z_N$
proves the theorem.
\end{proof}

\section{Final remarks}\label{sec:remarks}

We first record the following bonus from the orientation double
covers of these examples.

\begin{addendum}\label{add:orientation-cover}
For every positive integer $N$, there exists a closed orientable
aspherical topological $4$--manifold $\widetilde Z_N$ of signature
zero which admits at least $N$ pairwise nondiffeomorphic smooth
structures.
\end{addendum}

\begin{proof}
Fix $N$ and take the orientation double covers $\widetilde Y_i$
of the $Y_i$ constructed in the proof of Theorem~\ref{thm:main}.

Since $\rho(H)=\{0,a\}$ and $a$ has odd weight, $C$ is exactly the
subgroup of $H$ that preserves the orientation of $D_i$. Thus
\[
\widetilde Y_i=D_i/C\longrightarrow D_i/H=Y_i
\]
is the orientation double cover, with the chamber pairing shown in
Figure~\ref{fig:involution}. Its free deck involution reverses
orientation, so $\sigma(\widetilde Y_i)=0$. The closed orientable
aspherical manifolds $\widetilde Y_i$ are homeomorphic to a common
topological manifold $\widetilde Z_N$ and are pairwise nondiffeomorphic
by the same chamberwise homeomorphisms and lifting argument,
also applying~\cite[Lemma~5.3]{DHHRS-2026} with $C$ in place of $H$.
The genus argument also rules out orientation reversing
diffeomorphisms.
\end{proof}

The examples of Theorem~\ref{thm:main} and
Addendum~\ref{add:orientation-cover} also provide further counterexamples
to the smooth analogue of the Borel conjecture~\cite{DHHRS-2026}.

\smallskip
\begin{rmk}[Why switch ingredients?] \label{rmk:chamber-choice}
Our choice of quotient subgroups also works if we use the
Hayden--Piccirillo chambers~\cite{HaydenPiccirillo2025} as in
DHHRS~\cite{DHHRS-2026}, giving an exotic pair of closed orientable
aspherical $4$--manifolds with free orientation reversing involutions
and exotic nonorientable quotients. However, it is unclear whether
their variants remain homeomorphic as the number $m$ of added
positive clasps increases; cf.~\cite[Remarks~2.2 and~5.4(b)]{DHHRS-2026}.
We use Akbulut--Yasui's construction~\cite{AkbulutYasui2013} to obtain
$N$ homeomorphic chambers with distinct values of $\mu$ for every
prescribed positive integer $N$.
\end{rmk}

\begin{samepage}
\begin{rmk}[Irreducibility] \label{rmk:irreducibility}
A smooth $4$--manifold is \emph{irreducible} if every
smooth connected sum decomposition has a homotopy $4$--sphere as one
of its summands. Our examples are
irreducible by asphericity. Indeed, the separating $3$--sphere of a
connected sum decomposition is nullhomotopic and bounds a contractible
piece; see e.g. the argument in~\cite[proof of Proposition~6.1]{BaykurStipsiczSzabo2024}.
Capping this piece with a $4$--ball gives a homotopy $4$--sphere as one of the summands. See also~\cite[Theorem~5.1]{KonnoLinMukherjeeMunozEchaniz2026}.
\end{rmk}
\end{samepage}

\vspace{-0.5\baselineskip}
Our work leaves the following question open
(cf.~\cite[Problem~4.2]{K3}):

\vspace{-0.5\baselineskip}
\begin{question}\label{qst:infinite}
Does there exist a closed nonorientable topological $4$--manifold
that admits infinitely many pairwise nondiffeomorphic smooth structures?
\end{question}

I originally hoped to use the DHHRS approach to obtain an infinite
family of nonorientable exotica.
The difficulty with the genus
argument we adopted here is that a fixed chamber type puts a uniform upper bound on
$\mu$. To see this, fix a chamber $X$. Since $X$ embeds in $D^4$,
its intersection form is zero. Moreover, $X\simeq T^2$, so $H_1(X)$
is torsion-free. The exact sequence of the pair then gives
a surjection $H_2(\partial X)\to H_2(X)$. Thus we can choose a
connected oriented surface $\Sigma\subset\partial X$ of genus $g$
whose class has nonzero image in $H_2(X)$. Now let $X'$ be any smoothing of the same topological chamber.
The restriction of a homeomorphism $X\to X'$ to the boundary is
isotopic to a diffeomorphism. We can therefore carry $\Sigma$ to
$\partial X'$ and push it into the interior, obtaining a smooth
surface of genus $g$ representing a nonzero class. Hence, $\mu(X')\le g$
for every such smoothing and only finitely many integer values of $\mu$
are possible for this fixed chamber type. Since our proof distinguishes
the chambers by strictly increasing values of $\mu$, we conclude that their common
homeomorphism type cannot stay fixed as $N$ grows arbitrarily large---which limits the genus argument.

\begin{rmk}[Some ideas that don't work] \label{rmk:earlier-methods}
I imagine that, like me, the reader may wonder why the methods that
already gave exotic pairs of closed nonorientable $4$--manifolds have
not also produced infinite families. I cannot give a full account of
this long and rich history here, but here are a few notes to myself
that I can share with the reader.

The $\eta$-invariant~\cite{Gilkey1985,Stolz1988}, for instance, distinguishes
several classical exotic pairs. These include the exotic projective
space of Fintushel--Stern, originally distinguished using a Rohlin-type
invariant and an equivariant signature defect~\cite[Proposition~1]{FintushelStern1981};
see~\cite[Theorem~A]{Oledzki1991}. It is additive modulo $2\Z$ under connected sums, for
chosen $\mathrm{Pin}^{+}$ structures on the
summands~\cite[Proposition~3]{Torres2017}. If replacing each summand
by its exotic copy changes $\eta$ by $1$, the invariant of the sum
therefore only records the parity of the number of replacements.

The Cappell--Shaneson invariant, defined in certain cases by a
signature-corrected Rohlin invariant, is unchanged by connected sum
with $S^2\x S^2$~\cite{RubermanSaveliev2007}.
More generally, invariants unchanged by taking connected sums with
$S^2\x S^2$ (e.g. the bordism invariants
in~\cite{Kreck1984,Debray2022}) cannot distinguish more than two smooth structures on a
fixed closed nonorientable $4$--manifold. Indeed, within a fixed
homeomorphism type there are at most two stable diffeomorphism
classes when the universal cover is spin, and only one when it is
nonspin~\cite[Remark~1.1, Theorem~3.2 and the proof of Proposition~4.1]{KasprowskiPowell2025}. This also covers the normal bordism invariants used to
detect examples admitting neither $\mathrm{Pin}^+$ nor $\mathrm{Pin}^-$
structures~\cite[Remark~1.3]{KasprowskiPowell2025}.
The smooth structures we construct on each $Z_N$ lie in a single
stable diffeomorphism class, so the classical invariants discussed
above cannot distinguish them. This differs from the classical exotic
pairs of Cappell--Shaneson, Fintushel--Stern, Akbulut and Kreck,
and the more recent examples of Torres and Bais--Torres, which remain
nondiffeomorphic after arbitrarily many $S^2\x S^2$
stabilizations; see~\cite{FintushelStern1981,KasprowskiPowell2025,TorresPrescribed2026,BaisTorresCutPaste2025}.
They therefore remain nondiffeomorphic after taking a product with $\R$.

A naive idea for obtaining infinite families of nonorientable
exotica is to take connected sums $Y\#M_i$, where $Y$ is a fixed
closed nonorientable $4$--manifold and the $M_i$ form an infinite
family of pairwise exotic closed orientable $4$--manifolds. One would then try to distinguish the smooth structures through their orientation double covers. The orientation double cover is $\widetilde Y\#M_i\#(-M_i)$, where
$\widetilde Y$ is the orientation double cover of $Y$.
Its ordinary Seiberg--Witten invariants vanish when
$b_2^+(M_i),b_2^-(M_i)>0$, and under presumably mild additional assumptions, its
Bauer--Furuta invariants vanish as
well~\cite[Lemma~2.6 and Theorem~2.8]{Yasui2019}.

There are possibly other ideas we are not covering here; also see~\cite[Question~D]{Torres2026} for a proposed family.
\end{rmk}

\newpage
\begin{samepage}
\subsection*{AI disclosure}
OpenAI's ChatGPT was used for figure generation, bibliography compilation and proofreading. It was also used for what does not appear in this article: before posting it on the arxiv preprint server, the author asked ChatGPT whether using the full genus functions of the chambers could yield an infinite exotic family, and it did no better than the author. The author takes full responsibility for all the writing, mathematical arguments and results.

\subsection*{Acknowledgements}
The author was partially supported by the NSF grant DMS-2506431
and by a Simons Foundation Travel Grant. We thank Matthias Kreck, Tye Lidman, Danny Ruberman, Nathan Sunukjian, and Kouichi Yasui
for helpful conversations and comments on a draft of this paper, Rafael Torres for pointing out his more recent preprints on the subject, and Bob Gompf for suggesting the stable equivalence conclusions.
\end{samepage}

\bibliographystyle{amsplain}

\bibliography{references}

@article{AkbulutYasui2013,
  author = {Akbulut, Selman and Yasui, Kouichi},
  title = {Cork twisting exotic {Stein} 4-manifolds},
  journal = {J. Differential Geom.},
  volume = {93},
  number = {1},
  pages = {1--36},
  year = {2013}
}

@article{DHHRS-2026,
  author = {Davis, Michael and Hayden, Kyle and Huang, Jingyin and Ruberman, Daniel and Sunukjian, Nathan},
  title = {Exotic aspherical 4-manifolds},
  journal = {Ann. of Math. (2)},
  volume = {204},
  number = {2},
  pages = {777--794},
  year = {2026},
  doi = {10.4007/annals.2026.204.2.6}
}

@unpublished{Torres2026,
  author = {Torres, Rafael},
  title = {A novel construction of homotopy 4-spheres and real projective 4-spaces {\`a} la {Cappell--Shaneson}},
  year = {2026},
  note = {arXiv:2607.15062}
}

@unpublished{BaisTorresCutPaste2025,
  author = {Bais, Valentina and Torres, Rafael},
  title = {A cut-and-paste mechanism to introduce fundamental group and construct new four-manifolds},
  year = {2025},
  note = {\href{https://arxiv.org/abs/2505.23167}{arXiv:2505.23167}}
}

@unpublished{TorresClassification2026,
  author = {Torres, Rafael},
  title = {Topological classification of certain nonorientable 4-manifolds with cyclic fundamental group of order {2 mod 4}},
  year = {2026},
  note = {\href{https://arxiv.org/abs/2601.05132}{arXiv:2601.05132}}
}

@article{CappellShaneson1976,
  author = {Cappell, Sylvain E. and Shaneson, Julius L.},
  title = {Some new four-manifolds},
  journal = {Ann. of Math. (2)},
  volume = {104},
  number = {1},
  pages = {61--72},
  year = {1976}
}

@article{FintushelStern1981,
  author = {Fintushel, Ronald and Stern, Ronald J.},
  title = {An exotic free involution on {$S^4$}},
  journal = {Ann. of Math. (2)},
  volume = {113},
  number = {2},
  pages = {357--365},
  year = {1981}
}

@incollection{Akbulut1984,
  author = {Akbulut, Selman},
  title = {A fake 4-manifold},
  booktitle = {Four-manifold theory ({Durham}, {N.H.}, 1982)},
  series = {Contemp. Math.},
  volume = {35},
  publisher = {Amer. Math. Soc.},
  address = {Providence, RI},
  pages = {75--141},
  year = {1984}
}

@article{Akbulut1988,
  author = {Akbulut, Selman},
  title = {Constructing a fake 4-manifold by {Gluck} construction to a standard 4-manifold},
  journal = {Topology},
  volume = {27},
  number = {2},
  pages = {239--243},
  year = {1988}
}

@incollection{FintushelStern1984,
  author = {Fintushel, Ronald and Stern, Ronald J.},
  title = {Another construction of an exotic {$S^1\mathbin{\widetilde{\times}}S^3\#S^2\times S^2$}},
  booktitle = {Four-manifold theory ({Durham}, {N.H.}, 1982)},
  series = {Contemp. Math.},
  volume = {35},
  publisher = {Amer. Math. Soc.},
  address = {Providence, RI},
  pages = {269--275},
  year = {1984}
}

@incollection{Kreck1984,
  author = {Kreck, Matthias},
  title = {Some closed 4-manifolds with exotic differentiable structure},
  booktitle = {Algebraic topology, {Aarhus} 1982},
  series = {Lecture Notes in Math.},
  volume = {1051},
  publisher = {Springer},
  address = {Berlin},
  pages = {246--262},
  year = {1984}
}

@article{Gompf1984,
  author = {Gompf, Robert E.},
  title = {Stable diffeomorphism of compact 4-manifolds},
  journal = {Topology Appl.},
  volume = {18},
  number = {2--3},
  pages = {115--120},
  year = {1984}
}

@article{Gompf1991,
  author = {Gompf, Robert E.},
  title = {Killing the {Akbulut--Kirby} 4-sphere, with relevance to the {Andrews--Curtis} and {Schoenflies} problems},
  journal = {Topology},
  volume = {30},
  number = {1},
  pages = {97--115},
  year = {1991}
}

@article{Gilkey1985,
  author = {Gilkey, Peter B.},
  title = {The eta invariant for even-dimensional {$\mathrm{Pin}^c$} manifolds},
  journal = {Adv. in Math.},
  volume = {58},
  number = {3},
  pages = {243--284},
  year = {1985}
}

@article{Stolz1988,
  author = {Stolz, Stephan},
  title = {Exotic structures on 4-manifolds detected by spectral invariants},
  journal = {Invent. Math.},
  volume = {94},
  number = {1},
  pages = {147--162},
  year = {1988}
}

@article{Oledzki1991,
  author = {Oledzki, Wies{\l}aw J.},
  title = {On the eta-invariant of {$\mathrm{Pin}^+$}-operator on some exotic 4-dimensional projective space},
  journal = {Compositio Math.},
  volume = {78},
  number = {1},
  pages = {1--27},
  year = {1991}
}

@article{HambletonKreckTeichner1994,
  author = {Hambleton, Ian and Kreck, Matthias and Teichner, Peter},
  title = {Nonorientable 4-manifolds with fundamental group of order 2},
  journal = {Trans. Amer. Math. Soc.},
  volume = {344},
  number = {2},
  pages = {649--665},
  year = {1994}
}

@article{Wang1995,
  author = {Wang, Zhenghan},
  title = {Classification of closed nonorientable 4-manifolds with infinite cyclic fundamental group},
  journal = {Math. Res. Lett.},
  volume = {2},
  number = {3},
  pages = {339--344},
  year = {1995}
}

@article{Torres2017,
  author = {Torres, Rafael},
  title = {Smooth structures on nonorientable four-manifolds and free involutions},
  journal = {J. Knot Theory Ramifications},
  volume = {26},
  number = {13},
  pages = {1750085},
  year = {2017}
}

@article{BaisTorres2025,
  author = {Bais, Valentina and Torres, Rafael},
  title = {Smooth structures on nonorientable 4-manifolds via twisting operations},
  journal = {Bull. Lond. Math. Soc.},
  volume = {57},
  number = {6},
  pages = {1768--1790},
  year = {2025},
  doi = {10.1112/blms.70060}
}

@article{TorresPrescribed2026,
  author = {Torres, Rafael},
  title = {Exotic non-orientable four-manifolds with prescribed fundamental group},
  journal = {Math. Proc. Cambridge Philos. Soc.},
  volume = {181},
  number = {1},
  pages = {773--780},
  year = {2026},
  doi = {10.1017/S0305004126101935}
}

@unpublished{KasprowskiPowell2025,
  author = {Kasprowski, Daniel and Powell, Mark},
  title = {Stably exotic 4-manifolds},
  year = {2025},
  note = {arXiv:2508.10499v2}
}

@article{Yasui2019,
  author = {Yasui, Kouichi},
  title = {Geometrically simply connected 4-manifolds and stable cohomotopy {Seiberg--Witten} invariants},
  journal = {Geom. Topol.},
  volume = {23},
  number = {5},
  pages = {2685--2697},
  year = {2019},
  doi = {10.2140/gt.2019.23.2685}
}

@article{Davis1983,
  author = {Davis, Michael W.},
  title = {Groups generated by reflections and aspherical manifolds not covered by {Euclidean} space},
  journal = {Ann. of Math. (2)},
  volume = {117},
  number = {2},
  pages = {293--324},
  year = {1983}
}

@article{HaydenPiccirillo2025,
  author = {Hayden, Kyle and Piccirillo, Lisa},
  title = {The trace embedding lemma and spinelessness},
  journal = {J. Differential Geom.},
  volume = {131},
  number = {1},
  pages = {133--163},
  year = {2025}
}

@article{FintushelStern1998,
  author = {Fintushel, Ronald and Stern, Ronald J.},
  title = {Knots, links, and 4-manifolds},
  journal = {Invent. Math.},
  volume = {134},
  number = {2},
  pages = {363--400},
  year = {1998},
  doi = {10.1007/s002220050268}
}

@article{Wall1967,
  author = {Wall, C. T. C.},
  title = {Locally flat {PL} submanifolds with codimension two},
  journal = {Proc. Cambridge Philos. Soc.},
  volume = {63},
  number = {1},
  pages = {5--8},
  year = {1967}
}

@article{Debray2022,
  author = {Debray, Arun},
  title = {Stable diffeomorphism classification of some unorientable 4-manifolds},
  journal = {Bull. Lond. Math. Soc.},
  volume = {54},
  number = {6},
  pages = {2219--2231},
  year = {2022},
  doi = {10.1112/blms.12688}
}

@article{Park2002,
  author = {Park, Jongil},
  title = {The geography of spin symplectic 4-manifolds},
  journal = {Math. Z.},
  volume = {240},
  number = {2},
  pages = {405--421},
  year = {2002}
}

@article{AkhmedovParkSakalli2023,
  author = {Akhmedov, Anar and Park, B. Doug and Sakall{\i}, S{\"u}meyra},
  title = {Exotic smooth structures on connected sums of {$S^2\times S^2$}},
  journal = {Q. J. Math.},
  volume = {74},
  number = {2},
  pages = {799--809},
  year = {2023},
  doi = {10.1093/qmath/haad002}
}

@article{BaykurHamada2026,
  author = {Baykur, R. {\.{I}}nan{\c{c}} and Hamada, Noriyuki},
  title = {Exotic 4-manifolds with signature zero},
  journal = {Selecta Math. (N.S.)},
  volume = {32},
  pages = {Paper No. 74},
  year = {2026},
  doi = {10.1007/s00029-026-01177-3}
}

@unpublished{BaykurStipsiczSzabo2024,
  author = {Baykur, R. {\.{I}}nan{\c{c}} and Stipsicz, Andr{\'a}s I. and Szab{\'o}, Zolt{\'a}n},
  title = {Smooth structures on four-manifolds with finite cyclic fundamental groups},
  year = {2024},
  note = {arXiv:2406.09007}
}

@incollection{FintushelStern2009,
  author = {Fintushel, Ronald and Stern, Ronald J.},
  title = {Six lectures on four 4-manifolds},
  booktitle = {Low dimensional topology},
  series = {IAS/Park City Math. Ser.},
  volume = {15},
  publisher = {Amer. Math. Soc.},
  address = {Providence, RI},
  pages = {265--315},
  year = {2009},
  doi = {10.1090/pcms/015/09}
}

@article{BaykurHamada2024,
  author = {Baykur, R. {\.{I}}nan{\c{c}} and Hamada, Noriyuki},
  title = {Lefschetz fibrations with arbitrary signature},
  journal = {J. Eur. Math. Soc. (JEMS)},
  volume = {26},
  number = {8},
  pages = {2837--2895},
  year = {2024},
  doi = {10.4171/JEMS/1326}
}

@article{PrzytyckiSwiatkowski2009,
  author = {Przytycki, Piotr and {\'S}wi{\k{a}}tkowski, Jacek},
  title = {Flag-no-square triangulations and {Gromov} boundaries in dimension 3},
  journal = {Groups Geom. Dyn.},
  volume = {3},
  number = {3},
  pages = {453--468},
  year = {2009},
  doi = {10.4171/GGD/66}
}

@book{OzbagciStipsicz2004,
  author = {Ozbagci, Burak and Stipsicz, Andr{\'a}s I.},
  title = {Surgery on contact 3-manifolds and {Stein} surfaces},
  series = {Bolyai Soc. Math. Stud.},
  volume = {13},
  publisher = {Springer},
  address = {Berlin},
  year = {2004},
  doi = {10.1007/978-3-662-10167-4}
}

@unpublished{YasuiGenus2025,
  author = {Yasui, Kouichi},
  title = {Corks, exotic 4-manifolds and genus functions},
  year = {2025},
  note = {arXiv:2501.18584v5}
}

@book{Davis2025,
  author = {Davis, Michael W.},
  title = {The Geometry and Topology of {Coxeter} Groups},
  edition = {2nd},
  series = {Springer Monographs in Mathematics},
  publisher = {Springer},
  address = {Cham},
  year = {2025},
  doi = {10.1007/978-3-031-91303-7}
}

@article{RubermanSaveliev2007,
  author = {Ruberman, Daniel and Saveliev, Nikolai},
  title = {Dirac operators on manifolds with periodic ends},
  journal = {J. G\"okova Geom. Topol.},
  volume = {1},
  pages = {33--50},
  year = {2007}
}

@unpublished{KonnoLinMukherjeeMunozEchaniz2026,
  author = {Konno, Hokuto and Lin, Jianfeng and Mukherjee, Anubhav and Mu{\~n}oz-Ech{\'a}niz, Juan},
  title = {Exotic diffeomorphisms on closed aspherical {$4$}-manifolds},
  year = {2026},
  note = {\href{https://arxiv.org/abs/2610.21083}{arXiv:2610.21083}}
}

@book{K3,
  editor    = {Baykur, R. {\.{I}nan\c{c}} and Kirby, Robion C.
               and Ruberman, Daniel},
  title     = {{K3}: A New Problem List in Low-Dimensional Topology},
  series    = {Mathematical Surveys and Monographs},
  volume    = {295},
  publisher = {American Mathematical Society},
  address   = {Providence, RI},
  year      = {2026},
  note      = {Preliminary version available at
               \url{https://bpb-us-e2.wpmucdn.com/websites.umass.edu/dist/b/22144/files/2026/04/K3-problem-list-watermarked.pdf}}
}

@article{Ruberman1984,
  author = {Ruberman, Daniel},
  title = {Invariant knots of free involutions of {$S^4$}},
  journal = {Topology Appl.},
  volume = {18},
  number = {2--3},
  pages = {217--224},
  year = {1984},
  doi = {10.1016/0166-8641(84)90011-7}
}

@article{HausmannJahren2018,
  author = {Hausmann, Jean-Claude and Jahren, Bj{\o}rn},
  title = {A simplification problem in manifold theory},
  journal = {Enseign. Math.},
  volume = {64},
  number = {1--2},
  pages = {207--248},
  year = {2018},
  doi = {10.4171/LEM/64-1/2-8}
}

\end{document}